%Paper: alg-geom/9401004
%From: krasinsk%PLUNLO51.BITNET@SEARN.SUNET.SE
%Date: Fri, 21 Jan 1994 11:18:19 CET

                        %%%% PAPER %%%%

%AmS-TeX, Version 2.0,  9 pages

\input amstex
       \documentstyle{amsppt}
       \hoffset=.6in
       \NoBlackBoxes

\define\m{m\!\!-\!\!1}
\define\mj{m\!\!-\!\!1}
\define\md{m\!\!-\!\!2}

       \define\lo{\longrightarrow}
       
       \define\cv{\operatorname{const.}}

       \define\vct#1#2{#1_0,\ldots,#1_#2}
       \define\cuc{\Bbb C}
       \define\cac{\Bbb C^2}
       \define\cxy{\Bbb C[x,y]}

       \define\lam{\lambda}

       \font\bbf=cmr10 scaled \magstep 2
       \topmatter
       \address
        Institute of Mathematics,
        University of \L \'od\'z,
        S. Banacha 22 ,
       90-238 \L \'od\'z, Poland
       \endaddress
       \email krasinsk\@plunlo51.bitnet
       \endemail
       \endtopmatter
       \document
       \centerline{\bbf MATRIX IDENTITIES CONNECTED}
\vskip .2in
\centerline{\bbf WITH THE JACOBIAN CONJECTURE}

      \vskip .5in
       \centerline{\bbf T. Krasi\'nski}
       \vskip 1 in
       \font\rm=cmr10 scaled \magstephalf
       \font\it=cmti10 scaled \magstephalf
       \rm{

       \head 1. Introduction
       \endhead
       Let $f \in \cxy$ be a polynomial in two variables with complex
       coefficients. $f$ is said to be a {\it component of an automorphism}
       if there exists a polynomial $g \in \cxy $ such that $F:=(f,g):\cac
       \lo \cac$ is a polynomial automorphism of $\cac$ (i.e. there exists
       $F^{-1}$ and it is also polynomial).

In turn, $f$ is said to be a {\it Keller's component} if there exists $g\in
\cxy$ such that the jacobian Jac$(f,g)$ of the mapping $(f,g):\cac \lo \cac$
is a non-zero constant. Then $g$ is called {\it associated with} $f$ and
the pair $(f,g)$ a {\it Keller's mapping}.

It is well-known that the famous (and unsolved so far) jacobian conjecture
can be formulated in the following way:

(1) {\it If $f$ is a Keller's component,  then $f$ is a component of an
automorphism.}

It is easy to show that it sufficies to prove (1) for $f\in\cxy$ satisfying the
 following three conditions which will be called in the sequel the main
assumptions:

1. $f$ is of the form
$$f(x,y)=y^m+a_1(x)y^{m-1}+\hdots +a_m(x),\  m\geqslant 2,\
\text{deg}a_i\leqslant i,\
 i=1,\hdots ,m,$$

2. for every $ \lam \in \cuc$, $f-\lam $ is reduced in $\cxy$ ,

3. $\deg _yf'_x >0$.

In the paper, using one of the main results in [3], we give conditions for  $f$
(which satisfy the main assumptions) equivalent
to the fact that $f$ is a component of an
automorphism and a Keller's component, respectively (in the latter case under
some additional assumptions on $f$). These conditions are differential matrix
identities connected with the matrix
$$
M=
\bmatrix
 a'_1 & a'_2& \hdots &a'_m&&\\
&\ddots&&&\ddots&\\
&&a'_1&a'_2&\hdots&a'_m\\
m&(\m )a_1&\hdots&a_{\m }&&\\
&\ddots&&&\ddots&\\
&&m&(\m )a_1&\hdots &a_{\m }
\endbmatrix
\aligned
       & { \left.           \aligned
      &  \\
      &  \endaligned \right\} \, (m-1)-\text{rows}}\\
  & { \left.      \aligned
& \\
      &  \endaligned \right\} \, (m-1)-\text{rows}},
\endaligned
$$
which determinant is equal, up to a non-zero constant,
 to the resultant $R_y(f'_x,f'_y) $ of the partial
derivativies $f'_x,f'_y$ of $f$ with respect to the variable $y$.
Precisely, if $a'_1=0,\hdots, a'_k=0,\ a'_{k+1}\ne 0$, then
$$
\text{det}M=(-1)^{k(m+1)}m^kR_y(f'_x,f'_y).
$$
In the sequel, for a square matrix $A=
[a_{ij}]_{1\leqslant i,j \leqslant n},\  a_{ij}
\in \cuc [x]$, we shall denote by $|a_{ij}|_{x^k}$ k-th derivative of the
determinant of $A$.
\head {2.Component of an automorphism}
\endhead
\proclaim{Theorem A} Let $f \in \cxy$ satisfy the main assumptions. Then $f$ is
a component of an automorphism if and only if the following identities hold
$$ \sum ^{m-1}\Sb
r_1,\hdots ,r_i = 1\\
s_1,\hdots,s_j=1\\
r_1<\hdots <r_i\\
s_1<\hdots <s_j\\
\endSb
\vmatrix
a'_1&a'_2&\hdots &a'_m&&&&&\\
&&\ddots &&&&&&\\
&0&0&\hdots& 1&&&&\\
&&&\ddots &&&&&\\
&&&a'_1&a'_2&\hdots&a'_m&&\\
&&&&&\ddots &&&\\
&&&&0&0&\hdots&1&\\
&&&&&&\ddots&&\\
&&&&&a'_1&a'_2&\hdots&a'_m\\
m&(\m)a_1&\hdots &a_{\m}&&&&&\\
&&\ddots &&&&&&\\
&0&0&\hdots& 1&&&&\\
&&&\ddots &&&&&\\
&&&m&(\m)a_1&\hdots&a_{\m}&&\\
&&&&&\ddots &&&\\
&&&&0&0&\hdots&1&\\
&&&&&&\ddots&&\\
&&&&&m&(\m)a_1&\hdots&a_{\m}
\endvmatrix
\matrix
 &\\
 &\\
& \\
\leftarrow&r_1\\
 &\\
\vdots&\\
 \\
\leftarrow&r_i\\
 &\\
 &\\
 &\\
 &\\
&\\
\leftarrow&(\m)+s_1\\
 &\\
\vdots&\\
 &\\
 &\\
\leftarrow&(\m )+s_j\\
 &\\
 x^k&\\
\endmatrix
\equiv 0
\tag{A}
$$
for $k=1,\hdots ,m-1,\quad  i,j=0,1,\hdots, k-1,\quad i+j=k-1 $
(the above determinant is obtained from the matrix $M$ by replacing
the rows with indices $r_1,\hdots ,r_i,(m-1)+s_1,\hdots,(m-1)+s_j$
by appropriate versors $[0,\hdots,0,1,0,\hdots,0]$ with the unit
at $m+r_1-1,\hdots,m+r_i-1,m+s_1-1,\hdots,m+s_j-1$ place, respectively).

\endproclaim
\remark{Remark 1} In particular, for $k=1$ we have the identity
$$
\vmatrix
 a'_1 & a'_2& \hdots &a'_m&&\\
&\ddots&&&\ddots&\\
&&a'_1&a'_2&\hdots&a'_m\\
m&(\m)a_1&\hdots&a_{\m}&&\\
&\ddots&&&\ddots&\\
&&m&(\m)a_1&\hdots &a_{\m}
\endvmatrix _x\equiv 0,$$
and for $k=2$ we have two identities
$$
\sum^{m-1}_{r=1}
\vmatrix
 a'_1 & a'_2& \hdots &a'_m&&&&\\
&\ddots&&&\ddots&&&\\
&&0&0&\hdots&1&&\\
&&&\ddots&&&\ddots&\\
&&&&a'_1&a'_2&\hdots&a'_m\\
m&(\m)a_1&\hdots&a_{\m}&&&&\\
&\ddots&&&\ddots&&&\\
&&\ddots&&&\ddots&&\\
&&&\ddots&&&\ddots&\\
&&&&m&(\m)a_1&\hdots &a_{\m}
\endvmatrix
\matrix & \\ &\\ & \\ \leftarrow&r\\ & \\ & \\ & \\&\\ & \\ & \\ & \\&
\\x^2&\endmatrix
\equiv 0,
$$
\vskip.3in
 $$
\sum^{m-1}_{s=1}
\vmatrix
 a'_1 & a'_2& \hdots &a'_m&&&&\\
&\ddots&&&\ddots&&&\\
&&\ddots&&&\ddots&&\\
&&&\ddots&&&\ddots&\\
&&&&a'_1&a'_2&\hdots&a'_m\\
m&(\m)a_1&\hdots&a_{\m}&&&&\\
&\ddots&&&\ddots&&&\\
&&0&0&\hdots&1&&\\
&&&\ddots&&&\ddots&\\
&&&&m&(\m)a_1&\hdots &a_{\m}
\endvmatrix
\matrix &\\ &\\ &\\&\\ & \\ & \\&\\&\\&\\ \leftarrow&(\m)+s
\\&\\&\\x^2&\endmatrix
\equiv 0,
$$
\endremark

\remark{Remark 2} Obviously in (A) we may additionally assume that
$r_p\ne s_q,\quad p=1,\hdots ,i,\quad q=1,\hdots ,j$.
\endremark
\remark{Remark 3} The identities (A) give also a set of equations describing
coefficients of components of polynomial automorphisms in the space of all
coefficients of a polynomial in two variables of degree $m$.
\endremark
\demo{Proof of Theorem A} Let $f$ satisfy the main assumptions.
Denote by $u,v$ new variables and consider the polynomials
$f'_x-u, f'_y-v \in \Bbb C[x,y,u,v]$. Let
$$
Q(x,u,v)= Q_N(u,v)x^N+Q_{N-1}(u,v)x^{N-1}+\hdots +Q_0(u,v)
$$
be the resultant of these polynomials with respect to the variable $y$.
By definition
$$
Q(x,u,v)=R_y(f'_x-u,f'_y-v)=C
\vmatrix
 a'_1 & a'_2& \hdots &a'_m-u&&\\
&\ddots&&&\ddots&\\
&&a'_1&a'_2&\hdots&a'_m-u\\
m&(\m )a_1&\hdots&a_{\m }-v&&\\
&\ddots&&&\ddots&\\
&&m&(\m )a_1&\hdots &a_{\m}-v
\endvmatrix,\tag2
$$
where $C$ is a non-zero constant.
It is easy to check that $Q(x,u,v) \not\equiv 0$. Moreover, by main assumptions
1. and 2. we have that $R_y(f'_x,f'_y)\not\equiv 0 $, i.e.
$$
Q(x,0,0)\not\equiv 0.\tag3
$$

In  paper \cite{3}, Th.3.4 (see also \cite{2}, Th.5.2, where the inequality
$\text{ord}_{(0)}Q_N\ne 0$ should be replaced by the equality
$\text{ord}_{(0)}Q
_N
= 0$) the following characterization, in terms of $Q$, of a component of an
automorphism was given:
\vskip.2in
{\sl If $f\in\cxy,\quad 0<\text{\rm {deg}}f-1= \text{\rm{deg}}_yf'_y$ and
$0<\text{\rm{deg}}_yf'_x$, then $f$ is a
component of an automorphism if and only if $\text{\rm ord}_{(0,0)}Q_0=0$ and
$\text{\rm ord}_{(0,0)}
Q_i\ge i$ for $i=1,\dots ,N$.}
\vskip.2in
Actually, Theorem 3.4 in \cite{3}  is slightly weaker. However, from Remark
 11.2 and Theorem 9.4 of that paper we easily obtain the above form of
the theorem.

{}From the main assumptions on $f$ it follows that $f$ satisfies the general
assumptions of the cited theorem. So, $f$ is a component of an automorphism
if and only if $\text{ord}_{(0,0)}Q_0=0$ and
$\text{ord}_{(0,0)}Q_i\ge i$ for $i=1,\dots ,N$. These conditions are
equivalent to the following
$$
Q_{0}(0,0)\ne 0,\ \frac{\partial ^{i+j}Q_k}{\partial u^i\partial v^j}(0,0)=0,
\ k=1,\dots,N,\ i,j=0,1,\dots,k-1,\ i+j\le k-1.\tag4
$$
Since for every $k=1,\dots,N$ and $i,j\in \Bbb N$ we have
$$
\frac{\partial ^{i+j+k}Q}{\partial u^i\partial v^j\partial x^k}(x,0,0)=
\frac{N!}{(N-k)!}\frac{\partial ^{i+j}Q_N}{\partial u^i\partial
v^j}(0,0)x^{N-k}
+\dots +k!\frac{\partial ^{i+j}Q_k}{\partial u^i\partial v^j}(0,0),
$$
then (4) (taking into account (3)) is equivalent to the following:
$$
\frac{\partial ^{i+j+k}Q}{\partial u^i\partial v^j\partial x^k}(x,0,0)\equiv 0,
\ k=1,\dots,N,\ i,j=0,1,\dots,k-1,\ i+j=k-1.\tag5
$$
{}From form (2) of $Q$ it follows that we have obtained (A) for $k=1,\dots,N$.
Obviously, if $N\le m-1$ then from (5) we obtain that the identities in (A)
are satisfied for $k=N+1,\dots,(m-1)$. In turn, if $m-1<N$ then directly from
identities (A) it follows that they also hold for $k=m,\dots,N$. Hence (5) is
al
so
satisfied for $k=m,\dots,N$.

In consequence, (5) is equivalent to (A) in every case.

\enddemo

\head 3. Keller's component
\endhead
\proclaim{Theorem B} Let $f\in\cxy$ satisfy the main assumptions.
Then $f$ is a Keller's component for which there exists an associated
component $g$ such that $\text{\rm deg}_yg\le m $ if and only if the following
identities hold
$$
\align
&\vmatrix
 a'_1 & a'_2& \hdots &a'_m&&&&\\
&\ddots&&&\ddots&&&\\
&&0&0&\hdots&\hdots&0&1\\
&&&\ddots&&&\ddots&\\
&&&&a'_1&a'_2&\hdots&a'_m\\
m&(\m)a_1&\hdots&a_{\m}&&&&\\
&\ddots&&&\ddots&&&\\
&&\ddots&&&\ddots&&\\
&&&\ddots&&&\ddots&\\
&&&&m&(\m)a_1&\hdots &a_{\m}
\endvmatrix
\matrix & \\ & \\ & \\  \leftarrow&k\\ & \\ & \\ & \\&\\ & \\ & \\ & \\&
\\x&\endmatrix
+\tag{B} \\
(m-k)
&\vmatrix
 a'_1 & a'_2& \hdots &a'_m&&&&\\
&\ddots&&&\ddots&&&\\
&&\ddots&&&\ddots&&\\
&&&\ddots&&&\ddots&\\
&&&&a'_1&a'_2&\hdots&a'_m\\
m&(\m)a_1&\hdots&a_{\m}&&&&\\
&\ddots&&&\ddots&&&\\
&&0&0&\hdots&\hdots&0&1\\
&&&\ddots&&&\ddots&\\
&&&&m&(\m)a_1&\hdots &a_{\m}
\endvmatrix
\matrix &\\ &\\ &\\ & \\& \\ &\\&\\&\\&\\ \leftarrow&(\m)+(k\!\!-\!\!1)
\\&\\&\\
&\endmatrix
\equiv 0
\endalign
$$
for $k=0,\dots,m-1$, where for $k=0,1$ the above identities mean
$$
%\align
\vmatrix
 a'_1 & a'_2& \hdots &a'_m&&\\
&\ddots&&&\ddots&\\
&&a'_1&a'_2&\hdots&a'_m\\
m&(\m)a_1&\hdots&a_{\m}&&\\
&\ddots&&&\ddots&\\
&&m&(\m)a_1&\hdots &a_{\m}
\endvmatrix _x\equiv 0,
% \\
$$
$$
\vmatrix
0&0&\hdots&\hdots&\hdots&0&1\\
 0&a'_1 & a'_2& \hdots &a'_m&&&\\
&&\ddots&&&\ddots&\\
&&&a'_1&a'_2&\hdots&a'_m\\
m&(\m)a_1&\hdots&a_{\m}&&&&\\
&\ddots&&&\ddots&&&\\
&&\ddots&&&\ddots&&\\
&&&m&(\m)a_1&\hdots &a_{\m}
\endvmatrix _x\equiv 0.
%\endalign
$$
\endproclaim

\demo{Proof} 1. Let $f$ be a Keller's component for which
 there exists an associated component $g$ such that $\text{deg}_yg\le m$.
Then $g(x,y)=b_0(x)y^m+b_1(x)y^{m-1}+\dots+b_m(x),\ b_i\in \Bbb C[x],\
i=0,\dots,m$. It is easy to show that we can assume that $b_0\equiv 0$.
In fact, it is easy to see that from the jacobian condition $\text{Jac}
(f,g)=\cv \ne 0$, it follows that $b'_0(x)\equiv 0$. Hence $b_0\equiv 0
$ or $b_0\equiv
\cv \ne 0.$ In the first case the assertion is satisfied, whereas in the
 second it sufficies to replace $g$ by $g-b_0f$.

So, in the sequel we shall assume that
$$
g(x,y)=b_1(x)y^{m-1}+\dots+b_m(x),\ b_i\in \Bbb C[x],\
i=1,\dots,m.
$$
Additionally, we shall assume that Jac$(f,g)=1$. Then
$$
\align
1=f'_xg'_y-f'_yg'_x&=(a'_1y^{m-1}+a'_2y^{m-2}+\dots+a'_m)((m-1)b_1y^{m-2}
+\dots+b_{m-1}) \\
&-(my^{m-1}+(m-1)a_1y^{m-2}+\dots+a_{m-1})(b'_1y^{m-1}+\dots+b'_m).
\endalign
$$

Comparing the coefficients of the both sides of this equality we obtain the
following system of polynomial equations
       $$
       \xxalignat {7}
       0@!&=&&&&&&&&&&&-@!m@!&b'_1 \\
       0@!&=@!&(@!\mj@!)@!a'_{1}@!&b_1&&&&&&&
       -@!
m@!&b'_2&@!-@!(@!\mj@!)@!a_{1}&@!b'_1 \\
       0@!&=@!&(@!\mj@!)@!a'_{2}@!&b_1&@!+@!(@!\md@!)@!a'_{1}@!&b_2& &&&&
       -@!(@!\mj@!)@!a_{1}@!&b'_2&@!-@!(@!\md@!)@!a_{2}&@!b'_1 \\
       \hdotsfor{14} \\
       0@!&=@!&(@!\mj@!)@!a'_{\mj}@!&b_1&@!+@!(@!\md@!)@!a'_{\md}@!&b_2&@!+
       @!\dots@!+@!a'_1@!&b_{\mj}@!&-@!m@!&b'_m@!&-@!\dots@!-
       @!2@!a_{\md}@!&b'_2&@!-@!a_{\mj}&@!b'_1 \\
       0@!&=@!&(@!\mj@!)@!a'_{m}@!&b_1&@!+@!(@!\md@!)@!a'_{\mj}@!&b_2&@!+
       @!\dots@!+@!a'_2@!&b_{\mj}@!&-@!(@!\mj@!)@!a_1@!&b'_m@!&-@!\dots@!-
       @!a_{\mj}@!&b'_2& & \\
       \hdotsfor{14} \\
       1@!&=&&&&&a'_m@!&b_{\mj}@!&-@!a_{\mj}&b'_m&&&&
       \endxxalignat
       $$
Considering $b_i$, $b'_j,\ i=1,\dots,m-1,\ j=1,\dots,m$ as unknows we have
a square system of $2m-1$ equations. Simple calculations give that for the
 determinant $D$ of the system we have
$$
\align
D&=(-1)^{\frac{m(m-1)}{2}-1}m!
\vmatrix
 a'_1 & a'_2& \hdots &a'_m&&\\
&\ddots&&&\ddots&\\
&&a'_1&a'_2&\hdots&a'_m\\
m&(\m)a_1&\hdots&a_{\m}&&\\
&\ddots&&&\ddots&\\
&&m&(\m)a_1&\hdots &a_{\m}
\endvmatrix \\
&=(-1)^{\frac{m(m-1)}{2}-1}m!CR_y(f'_x,f'_y),
\endalign
$$
where $C$ is a non-zero constant.
Since $f$ is a Keller's component $R_y(f'_x,f'_y)=\cv\ne 0$. In
consequence, $D$ is a non-zero constant equal to $
(-1)^{\frac{m(m-1)}{2}-1}m!R$, where $R=CR_y(f'_x,f'_y)$. So, by the
 Cramer's rule we get formulae for solutions of this system.
Simple calculations give
$$
\align
b_i&=\frac{1}{(m-i)R}
\vmatrix
 a'_1 & a'_2& \hdots &a'_m&&&&\\
&\ddots&&&\ddots&&&\\
&&0&0&\hdots&\hdots&0&1\\
&&&\ddots&&&\ddots&\\
&&&&a'_1&a'_2&\hdots&a'_m\\
m&(\m)a_1&\hdots&a_{\m}&&&&\\
&\ddots&&&\ddots&&&\\
&&\ddots&&&\ddots&&\\
&&&\ddots&&&\ddots&\\
&&&&m&(\m)a_1&\hdots &a_{\m}
\endvmatrix
\matrix & \\ & \\ & \\  \leftarrow&i\\ & \\ & \\ & \\&\\ & \\ & \\ & \\& \\
&\endmatrix \\  i&=1,\dots,m-1,\tag7
\endalign
$$

$$
\align
b'_1&=0,\\
b'_j&=-\frac{1}{R}
\vmatrix
 a'_1 & a'_2& \hdots &a'_m&&&&\\
&\ddots&&&\ddots&&&\\
&&\ddots&&&\ddots&&\\
&&&\ddots&&&\ddots&\\
&&&&a'_1&a'_2&\hdots&a'_m\\
m&(\m)a_1&\hdots&a_{\m}&&&&\\
&\ddots&&&\ddots&&&\\
&&0&0&\hdots&\hdots&0&1\\
&&&\ddots&&&\ddots&\\
&&&&m&(\m)a_1&\hdots &a_{\m}
\endvmatrix
\matrix &\\ &\\ &\\ & \\& \\ &\\&\\&\\&\\
\leftarrow&(m\!\!-\!\!1)\!\!+\!\!(j\!\!-\!\!1) \\&\\&\\&\endmatrix, \\
j&=2,\dots,m.
\tag8\endalign
$$
{}From these formulae it follows that identities (B) are fulfilled for
$k=1,\dots, m-1$. Moreover, from the fact that $R_y(f'_x,f'_y)=
\cv\ne 0$ we obtain identity (B) for $k=0$.

2. Assume that identities (B) hold for a polynomial $f$. From the main
assumptions 1. and 2. we have that $R_y(f'_x,f'_y)\not\equiv 0$.
So, by the identity (B) for $k=0$ we obtain $R_y(f'_x,f'_y)=
\cv\ne 0$. Hence $\text{det}M=CR_y(f'_x,f'_y),\
C=\cv\ne 0$, is also a non-zero constant. Denote it by $R$.
Then we define polynomials $
b_i,\ i=1,\dots,m-1$ by formulae (7) and polynomials
$\widetilde b_j,\ j=1,\dots,m$ by formulae (8). Additionally, we put
$b_m=\int\widetilde b_m\,dx$ (an arbitrary indefinite integral of $\widetilde
 b_m)$. From identities (B) we have
$$
b'_i=\widetilde b_i,\quad i=1,\dots,m.\tag9
$$
Besides, from the formulae for $b_i$ and $\widetilde b_j,\ i=1,\dots,m-1,
\ j=1,\dots,m$, it follows that they are solutions of system (6).
 Hence and from (9) we obtain that for $g(x,y)=b_1(x)y^{m-1}+\dots+b_m(x) $
we have Jac$(f,g)=1$.

\enddemo

\remark{Remark 4} From the proof it follows that the assumption deg$a_i\le i,\
i=1,\dots,m$, in the theorem is superfluous.
\endremark
\remark{Remark 5} A characterization of polynomials, associated with a given
Keller's component, can be find in \cite1 .
\endremark
\head{4}
Corollaries
\endhead
\proclaim{Corollary 1} Let $l\in\Bbb N$. If for every $m\le l$ and
arbitrary polynomials $a_1,\dots,a_m\in\Bbb C[x]$ such that
{\rm deg}$a_i\le i,\
i=1,\dots,m$, identities (B) imply identities (A) then the jacobian
conjecture is true for polynomials of degree $\le l$.
\endproclaim
\demo{Proof} Let $(f,g)$ be an arbitrary Keller's mapping such that $m=
\text{deg}f \le l$ and $n=\text{deg}g\le l$. Using a linear change of
coordinates in $\cac$ we may assume that both $f$ and $g$ satisfy the main
assumptions. If $m\ge n$, then from Theorem B we get that the coefficients
 $a_1,\dots,a_m$ of $f$ satisfy (B). By assumption we obtain that (A)
 hold. Hence, by Theorem A, $f$ is a component of an automorphism. It
 easily implies that $(f,g)$ is a polynomial automorphism.

In the case $n>m$ considerations are similar.
\enddemo

{}From the above corollary we obtain directly
\proclaim{Corollary 2} If for arbitrary polynomials $a_1,\dots,a_m\in\Bbb C
[x]$ such that {\rm deg}$a_i\le i,\ i=1,\dots,m$, identities (B) imply
identitie
s
(A), then the jacobian conjecture is true.
\endproclaim

Before we pass to concrete cases let us notice
\remark{Remark 6} Without loss of generality we may add to the main assumptions
that
 $a_1\equiv 0$. It can be obtained by applying to $f$ the known triangular
automorphism $(x,y)\longmapsto (x,y-\frac1na_1(x))$ of $\cac$.
\endremark

Let us now consider two cases:\newline
\indent 1. $m=2$. It is easily seen then in this case identities (A) and (B)
are equivalent.
\newline
\indent 2. $m=3$. In this case, assuming that $a_1\equiv 0$ (see Remark 6)
identities (A) are as follows
$$
((a'_2)^2a_2+3(a'_3)^2)'\equiv 0,\quad a''_2\equiv 0,\quad a'''_3\equiv 0,
\tag10
$$
whereas identities (B) are
$$
((a'_2)^2a_2+3(a'_3)^2)'\equiv 0,\quad a''_2\equiv 0,\quad a''_3\equiv 0,
\tag11
$$
It is easy to show that conditions (10) and (11) are equivalent.

{}From these cases and Corollary 1 we obtain
\proclaim{Corollary 3} The jacobian conjecture is true for polynomials of
degree $\le 3$.
\endproclaim

\vskip.3in
\rightline{\it \L \'od\'z, December 1993}
\enddocument